\documentclass[aop,MSNbibl,seceqn,dvips]{arximspdf}
\usepackage{mathbh}

\doi{10.1214/11-AOP735} 
\volume{41}
\issue{3A}
\pubyear{2013}
\firstpage{1160}
\lastpage{1179}

\makeatletter

\newcommand{\rrvert}{\vert}
\newcommand{\llvert}{\vert}

\newtheorem{itlemma}{Lemma}[section]
\newtheorem{itproposition}[itlemma]{Proposition}

\newproclaim{itremark}[itlemma]{Remark}
\newproclaim{itremarks}[itlemma]{Remarks}
\newproclaim{itdefinition}[itlemma]{Definition}
\newproclaim{itexample}[itlemma]{Example}

\def\T{\mathcal{T}}
\def\A{\mathcal{A}}
\def\C{\mathcal{C}}
\def\S{\mathcal{S}}

\def\B{\mathbb{B}}
\def\N{\mathbb{N}}
\def\Z{\mathbb{Z}}
\def\R{\mathbb{R}}
\def\id{\mathbf{1}}

\def\p{\partial}
\def\bs{\setminus}

\def\ind{\mathbh{1}}
\makeatother

\begin{document}
\begin{frontmatter}

\title{Sublogarithmic fluctuations for internal DLA}
\runtitle{Fluctuations for internal DLA}

\begin{aug}
\author[A]{\fnms{Amine} \snm{Asselah}\corref{}\ead[label=e1]{amine.asselah@u-pec.fr}}
\and
\author[B]{\fnms{Alexandre} \snm{Gaudilli\`ere}\ead[label=e2]{gaudilli@cmi.univ-mrs.fr}}
\runauthor{A. Asselah and A. Gaudilli\`ere}
\affiliation{Universit\'e Paris-Est Cr\'eteil and Universit\'e de Provence}
\address[A]{LAMA\\
Universit\'e Paris-Est Cr\'eteil\\
61 avenue du g\'en\'eral de Gaulle\\
94010 Cr\'eteil cedex\\
France\\
\printead{e1}} 
\address[B]{LATP\\
Universit\'e de Provence\\
CNRS, 39 rue F. Joliot Curie\\
13013 Marseille\\
France\\
\printead{e2}}
\end{aug}

\received{\smonth{12} \syear{2010}}
\revised{\smonth{10} \syear{2011}}

%
\begin{abstract}
We consider internal diffusion limited aggregation in dimension
larger than or equal to two. This is a random cluster growth model,
where random walks start at the origin of the $d$-dimensional lattice,
one at a time,
and stop moving when reaching a site that is not occupied by previous walks.
It is known that the asymptotic shape of the cluster is
a sphere. When the dimension is two or more, we have shown in a
previous paper that
the inner (resp., outer) fluctuations of its radius
is at most of order $\log(\mathrm{radius})$ [resp., $\log^2(\mathrm
{radius})$].
Using the same approach, we improve the upper bound on the inner
fluctuation to
$\sqrt{\log(\mathrm{radius})}$ when $d$ is larger than or equal to three.
The inner fluctuation is then
used to obtain a similar upper bound on the outer fluctuation.
\end{abstract}

%

\end{frontmatter}

\section{Introduction}
This note is a companion to our paper \cite{AG}.
There, we introduced a family of cluster growth models with a
spherical asymptotic shape, but a wide
diversity of shape fluctuations. Internal diffusion limited aggregation
(internal DLA) was one member of this family. More precisely,
the internal DLA cluster of volume $N$, say $A(N)$,
is obtained inductively as follows. Initially, we assume
that the explored region is empty, that is, $A(0)=\varnothing$.
Then, consider $N$ independent discrete-time
random walks $S_1,\ldots,S_N$ starting from 0.
Assume $A(k-1)$ is obtained, and define
\begin{equation}
\label{time-settling} \qquad\tau_k=\inf \bigl\{ {t\ge0\dvtx S_k(t)
\notin A(k-1)} \bigr\}\quad \mbox{and}\quad A(k)=A(k-1)\cup\bigl\{S_k(
\tau_k)\bigr\}.
\end{equation}
We call explorers the random walks obeying the aggregation rule
(\ref{time-settling}).
We say that the $k$th explorer is \textit{settled} on $S_k(\tau_k)$ after time~$\tau_k$, and is \textit{unsettled} before time~$\tau_k$.
The cluster $A(N)$ is interpreted as
the positions of the $N$ settled explorers.

In this paper we show how the tools
developed in~\cite{AG} lead in dimension
\mbox{$d\geq3$} to sharper estimates
on the fluctuations of $A(N)$
with respect to its spherical asymptotic shape.\vadjust{\goodbreak}
We keep the notation of \cite{AG}, and recall the basic ones to make
the paper as self-contained as possible.
We denote with $\|\cdot\|$ the Euclidean norm on $\R^d$.
For any $x$ in $\R^d$ and $r$ in $\R$, set
\begin{equation}
\label{ball-dfn} B(x,r) = \bigl\{ y\in\R^d \dvtx  \|y-x\| < r \bigr\}
\quad\mbox{and}\quad \B(x,r) = B(x,r) \cap\Z^d.
\end{equation}
For $\Lambda\subset\Z^d$,
$|\Lambda|$ denotes the number of sites in $\Lambda$, and the
boundary of $\Lambda$ is $\p\Lambda=\{z\notin\Lambda\dvtx \exists
y\in\Lambda,
\|y-z\|=1\}$. For a simple random walk, let
$H(\Lambda)$ denotes its first hitting time of $\Lambda$.
The inner error $\delta_I(n)$ is such that
\begin{equation}
\label{def-inner} n-\delta_I(n)=\sup \bigl\{ {r \geq0\dvtx  \B(0,r)
\subset A\bigl(\bigl|\B (0,n)\bigr|\bigr)} \bigr\}.
\end{equation}
Also, the outer error $\delta_O(n)$ is such that
\begin{equation}
\label{def-outer} n+\delta_O(n)=\inf \bigl\{ {r \geq0\dvtx  A\bigl(\bigl|
\B(0,n)\bigr|\bigr)\subset\B (0,r)} \bigr\}.
\end{equation}
Our main result is as follows.
\begin{itproposition}\label{prop-opti}
There are constants $\{\alpha_d, \beta_d, d\ge3\}$
such that in dimension $d\ge3$,
with probability~1,
\begin{equation}
\label{theo-res2} \lim\sup\frac{\delta_I(n)}{\sqrt{\log(n)}}\le\alpha_d\quad \mbox{and}\quad
\lim\sup\frac{\delta_O(n)}{\sqrt{\log(n)}}\le\beta_d.
\end{equation}
\end{itproposition}
\begin{itremark}\label{rem-opti}
For $d=2$ we show, with similar computations, that there
are constants $\alpha_2,\beta_2$ such that, with probability~1,
\begin{equation}
\label{theo-res1} \lim\sup\frac{\delta_I(n)}{\log(n)}\le\alpha_2\quad \mbox{and}\quad
\lim\sup\frac{\delta_O(n)}{\log(n)}\le\beta_2.
\end{equation}
The inner error bound in~(\ref{theo-res1}) was already obtained in all
dimensions in~\cite{AG}.
Recently, Jerison, Levine and Sheffield~\cite{levine} established,
in dimension two and with a different method,
the estimates (\ref{theo-res1}).
Also, they announced in~\cite{levine} that the approach
they followed could be adapted in dimension $d\geq3$
to get~(\ref{theo-res2}).
\end{itremark}

Let us describe the main steps.
The inner error is at the heart of
the argument. It is based on a large deviation estimate which
refines our previous estimates, with interest of its own.
For a real $x$, let $\lfloor x\rfloor$ be the integer part of $x$.
\begin{itlemma}\label{lem-ld}
Choose $R$ and $A$ large enough. Assume that
$\lfloor A R^d\rfloor$ explorers lie initially on $\B(0,R/2)$. We
call $\eta$
the initial configuration of these explorers and $A(\eta)$
the cluster they produce. There are positive
constants $\{\kappa_d, d\ge2\}$ independent of $R$ and $A$, such
that when $d\ge3$,
\begin{equation}
\label{bound-d3} P \bigl( {\B(0,R)\not\subset A(\eta)} \bigr) \le\exp \bigl( {-
\kappa_d AR^2} \bigr),
\end{equation}
and when $d=2$, we have
\begin{equation}
\label{bound-d2} P \bigl( {\B(0,R)\not\subset A(\eta)} \bigr) \le\exp \biggl( {-
\kappa_2 \frac{A R^2}{\log(R)}} \biggr).\vadjust{\goodbreak}
\end{equation}
\end{itlemma}
\begin{itremark}\label{rem-ld}
The reason behind the previous lemma, in $d\ge3$, is that out of
$\lfloor A R^d\rfloor$ explorers, only about $AR^2$ eventually hit a
fixed site
on the boundary of $\B(0,R)$,
so that it is only
these very explorers that need to be pushed away from this very site. The cost
should be proportional to $AR^2$.
\end{itremark}

For the outer error, we use
a large deviation estimate symmetrical to Lemma~\ref{lem-ld}
as well as our coupling between internal DLA
and the \textit{flashing process} of~\cite{AG}.
The latter large deviation estimate
was recently proved by Jerison, Levine and Sheffield in~\cite{levine}.
\begin{itlemma}[(Lemma A of Jerison, Levine and Sheffield \cite{levine})]\label{lem-levine}
For $\beta$ and $R$ positive reals, assume that
$\lfloor\beta R^d\rfloor$ explorers lie initially outside $\B(0,R)$.
We call $\eta$
the initial configuration of these explorers and $A(\eta)$
the cluster they produce.
There are positive constants $\{\kappa_d', d\ge2\}$, such that
for $\beta$ small enough,
we have when $d\ge3$,
\begin{equation}
\label{levine-1} P \bigl( {0\in A(\eta)} \bigr)\le\exp\bigl(-
\kappa_d' R^2\bigr),
\end{equation}
whereas when $d=2$, we have
\begin{equation}
\label{levine-2} P \bigl( {0\in A(\eta)} \bigr)\le\exp \biggl( {-
\kappa_2' \frac
{R^2}{\log(R)}} \biggr).
\end{equation}
\end{itlemma}
We give an alternative proof of this result,
based on estimating
the probability of crossing a shell, while avoiding traps.
\begin{itlemma}\label{lem-ag} Consider $d\ge2$. Fix a positive real
$R$, and
start a random walk on $z\in\p\B(0,2R)$.
There are positive constants $\{\kappa_d,a_d\}$ such that for any
$V$ subset of the shell $\S=\B(0,2R)\bs\B(0,R)$, we have
\begin{equation}
\label{avoiding-traps} P_z \bigl( {H \bigl( {\B(0,R)} \bigr)<H
\bigl(V^c\bigr)} \bigr)\le\exp \biggl( {a_d-
\kappa_d \frac {R}{\rho}} \biggr)\qquad \mbox{where } \rho^{d-1}=
\frac{|V|}{R}.\hspace*{-35pt}
\end{equation}
\end{itlemma}
\begin{itremark}\label{rem-AG}
$V^c=\S\bs V$ is interpreted as traps. Note that $\rho$ is
proportional to
the radius of a cylinder of height $R$ and volume $|V|$.
We can also read (\ref{avoiding-traps}) in the following way:
\begin{equation}
\label{avoid-traps} P_z \bigl( {H\bigl(\B(0,R)\bigr)<H
\bigl(V^c\bigr)} \bigr)\le\exp \biggl( {a_d-
\kappa_d \biggl( {\frac{R^d}{|V|}} \biggr)^{{1}/{(d-1)}}} \biggr).
\end{equation}
This shows that for (\ref{avoid-traps}) to be an effective
inequality, one needs that $|V|$ be smaller than $R^d$.
The power $1/(d-1)$ on $R^d/|V|$ in (\ref{avoid-traps}) is not
important in proving Lemma~\ref{lem-levine}. If one were willing to
accept the weaker power $1/d$, then one would have the following simple
heuristics in dimension $d\ge3$.
Let $t$ denote the time the walk spends in the annulus of height $R$.
On one hand, the central limit scaling yields that\vadjust{\goodbreak}
this probability of such a stay is of order $\exp(-c R^2/t)$. On the
other hand,
all this time should be spent on sites of~$V$, and
it is well known that the probability
is of order $\exp(-\kappa_d t/|V|^{2/d})$. Putting together
these opposite requirements, and optimizing over $t$, we find
a statement weaker than (\ref{avoid-traps}), but sufficient for our present
purpose:
\begin{equation}
\label{traps-weak} P_z \bigl( {H\bigl(\B(0,R)\bigr)<H
\bigl(V^c\bigr)} \bigr)\le\exp \biggl( {a_d-
\kappa_d \biggl( {\frac{R^d}{|V|}} \biggr)^{{1}/{d}}} \biggr).
\end{equation}
Even though it is not written in \cite{AG},
inequality (\ref{traps-weak}) was the motivation behind the
introduction of flashing processes in \cite{AG}, which were
basically used to bypass this type of estimate.
In this paper we show how the use of flashing
explorers leads easily to Lemma~\ref{lem-ag}.
\end{itremark}
The rest of the paper is organized as follows.
In Section~\ref{sec-old}, we enounce some known results:
we recall the approach of Lawler, Bramson and Griffeath \cite{lawler92}
and useful large deviation estimates.
Then the inner error estimate is proved in Section~\ref{sec-inner}.
In Section~\ref{sec-outer}, we show how a flashing process permits
a simple control on the outer error. Finally, we have gathered in an
\hyperref[app]{Appendix} the proof of the large deviations Lemmas~\ref{lem-ld}, \ref
{lem-levine}
and~\ref{lem-ag}.

\section{Prerequisites}\label{sec-old}
\subsection{Notation}
We recall some notation of \cite{AG}.
The state space of configurations is $\N^{\Z^d}$,
its elements are denoted $\eta$ and they represent starting conditions
for a set of explorers, or random walks.
Two types of initial configurations play an important role here:
(i) the configuration $n\id_{z^*}$ formed by $n$ trajectories starting on
a given site $z^*$ and (ii) for $\Lambda\subset\Z^d$,
the configuration $\id_\Lambda$
that we simply identify with $\Lambda$.
For any configuration $\eta\in\N^{\Z^d}$, we write
\begin{equation}
\label{not-1} |\eta| = \sum_{z\in\Z^d} \eta(z).
\end{equation}
\begin{itdefinition}\label{def-W}
Let $R\in\R_+\cup\{\infty\}$.
For $z\in\B(0,R)\cup\p\B(0,R)$, we denote by
$M_R(\eta,z)$ [resp., $W_R(\eta,z)$]
the number of simple random walks
(resp., explorers) initially on $\eta$ that
hit $z$ when or before exiting $\B(0,R)$.
Thus, when $z\in\p\B(0,R)$,
$M_R(\eta,z)$ [resp., $W_R(\eta,z)$] is the number of
simple random walks (resp., explorers) which exit $\B(0,R)$
exactly on $z$.
\end{itdefinition}
\begin{itremark}
Note that trajectories of walkers and explorers
can be coupled to be the same up to the settling time
of the explorer, the walker then proceeding along its
simple random walk trajectories.
\end{itremark}
As in~\cite{lawler95} (Section~3),
it is useful to stop explorers
as they reach $\partial B(0,R)$,
for some $R>0$, and
then to define
$A_R(\eta)$ as the set of positions
of settled explorers.\vadjust{\goodbreak}
\begin{itdefinition}\label{def-M} Consider $R\in\R\cup\{\infty\}$.
We set
\begin{equation}
\label{def-Mtilde} \forall z\in\B(0,R)\qquad \tilde M_R(
\eta,z)=W_R(\eta,z)+M_R\bigl(A_R(\eta),z
\bigr).
\end{equation}
\end{itdefinition}
Finally, for any function $F\dvtx \Z^d\to\R$ and subset $\Lambda\subset\Z^d$, we
denote
\[
F(\Lambda)=\sum_{z\in\Lambda} F(z).
\]

\subsection{On a classical approach}
We recall the approach of Lawler, Bramson
and Griffeath in~\cite{lawler92}. Send $N=|\B(0,n)|$ explorers
from the origin.
The approach of \cite{lawler92} is based on the following observations.
(i) If explorers did not settle,
they would just be independent random walks; (ii) exactly one
explorer occupies each site of the cluster.
Then, observations (i) and (ii) imply that for any integer $n$
and $z\in\B(0,n)$,
\begin{equation}
\label{main-lawler} \tilde M_n(N\ind_0,z):=
W_n(N\ind_0,z)+M_n\bigl(A_n(N),z
\bigr)\stackrel{\mathrm{law}} {\ge}M_n(N\ind_0,z).
\end{equation}
When $z\in\p\B(0,n)$, inequality (\ref{main-lawler}) becomes an equality,
\begin{equation}
\label{main-border} W_n(N\ind_0,z)+M_n
\bigl(A_n(N),z\bigr)\stackrel{\mathrm{law}}{=}M_n(N
\ind_0,z).
\end{equation}
Note that for any set $\Lambda\subset\B(0,n)$, $M_n(\Lambda,z)$
is a sum of independent Bernoulli variables.
Note also that $A_n(N)\subset\B(0,n)$ so that for any
$z\in\B(0,n)\cup\p\B(0,n)$
\begin{equation}
\label{lower-lawler92} W_n(N\ind_0,z)+M_n
\bigl(\B(0,n),z\bigr)\ge\tilde M_n(N\ind_0,z).
\end{equation}
However, Lawler et al.
did not use that $W_n(N\ind_0,z)$ and $M_n(\B(0,n),z)$ were independent.
They could only obtain a rough estimate on the lower tail
of $W_n(N\ind_0,z)$. This in turn gave some estimates on the inner error,
which was used to derive bounds on the outer error,
by using
that the cluster covers $\B(0,n-\delta_I(n))$.
In other words, from (\ref{main-border}),
and the definition of $\delta_I(n)$, for $R>n$ and $z\in\p\B(0,R)$,
\begin{equation}
\label{upper-lawler92} W_R(N\ind_0,z)+M_R
\bigl(\B\bigl(0,n-\delta_I(n)\bigr),z\bigr)\le\tilde
M_R(N\ind_0,z).
\end{equation}
Therefore, if $\delta_I(n)$ is likely to be smaller than $r<n<R$, and
$z\in\p\B(0,R)$, we have
\begin{equation}
\label{upper-main} \ind_{\{\delta_I(n)\le r\}} \bigl( {W_R(N
\ind_0,z)+M_R\bigl(\B (0,n-r),z\bigr)\le\tilde
M_R(N\ind_0,z)} \bigr).
\end{equation}
We will also make use of the independence of
the $\sigma$-fields generated by the events
$\{\delta_I(n)\le r\}$
and the random variables $W_R(N\ind_0,z)$
on the one hand,
and that generated by the random variable
$M_R(\B(0,n-r),z)$ on the other.

\subsection{On sums of Bernoulli variables}
Let us now recall a simple tool of~\cite{AG} in estimating
deviations in view of (\ref{lower-lawler92}) and (\ref{upper-lawler92}).
We first enounce the lower tail estimate.\vadjust{\goodbreak}
\begin{itlemma}\label{lem-AG-lower}
Suppose that a sequence of
random variables $\{W_n,M_n,L_n,\break \tilde M_n, n\in\N\}$, and
a sequence of real numbers $\{c_n, n\in\N\}$,
satisfy for each $n\in\N$,
\begin{equation}
\label{tool-lower} W_n+L_n+c_n\ge\tilde
M_n\quad \mbox{and}\quad \tilde M_n\stackrel{\mathrm{law}}
{=}M_n.
\end{equation}
Assume that $W_n$ and $L_n$ are independent, and that
$L_n$ and $M_n$ both are sums of independent Bernoulli variables.
Assume that the
Bernoulli variables $\{Y^{(n)}_1,\ldots,Y^{(n)}_{N_n}\}$
whose sum is $L_n$, satisfy for some $\kappa>1$,
\begin{eqnarray*}
&&\mathrm{(H1)} \quad \sup_n \sup_{i\le N_n} E\bigl[Y^{(n)}_i
\bigr]<\frac{\kappa- 1}{\kappa},\\
&&\mathrm{(H2)} \quad\mu_n:=E[M_n]-E[L_n]
\geq0.
\end{eqnarray*}
Then,
for any $n$ in ${\mathbb N}$ and $\xi_n$ in ${\mathbb R}$,
we have for all $\lambda\geq0$,
\begin{equation}
\label{tool-2}\qquad P ( {W_n<\xi_n} )\le\exp \Biggl( {-
\lambda(\mu_n-\xi_n-c_n)+\frac {\lambda^2}{2}
\Biggl( {\mu_n+\kappa\sum_{i=1}^{N_n}
E\bigl[Y^{(n)}_i\bigr]^2} \Biggr)} \Biggr).
\end{equation}
\end{itlemma}
The upper tail estimate needs other assumptions.
\begin{itlemma}\label{lem-AG-upper}
Assume for each $n\in\N$, and for an event $\A_n$,
\begin{equation}
\label{tool-upper} \ind_{\A_n} ( {W_n+L_n} )\le
\tilde M_n \quad\mbox {and} \quad\tilde M_n\stackrel{
\mathrm{law}} {=}M_n.
\end{equation}
Assume that $W_n$ and $L_n$ are independent, $\ind_{\A_n}$ and $L_n$
are independent and that
$L_n$ and $M_n$ both are sums of independent Bernoulli variables
such that $\mu_n:=E[M_n] - E[L_n]\ge0$.
Then, for all $n$ in ${\mathbb N}$, $\xi_n$ in ${\mathbb R}$
and $\lambda\in[0,\log2]$,
\begin{equation}
\label{toll-up} \qquad P ( {W_n\ge\xi_n, \A_n} )
\le\exp \biggl( {- \lambda (\xi_n-\mu_n)+
\lambda^2 \biggl( {\mu_n+4\sum
_{i} E\bigl[Y^{(n)}_i
\bigr]^2} \biggr)} \biggr).
\end{equation}
\end{itlemma}
\begin{itremark}\label{rem-1} This lower (resp., upper) tail estimate
turns out to be useful when $\xi_n + c_n$ is less than
(resp., $\xi_n$ is more than) $E[M_n] - E[L_n]$.
By Lemmas~\ref{lem-AG-lower} and~\ref{lem-AG-upper}
tail estimates reduce to a three-step strategy:
(i) estimation of $E[M_n] - E[L_n]$;
(ii) estimation of $\sum_i E^2[Y_i^{(n)}]$;
(iii) optimization in $\lambda$.
We emphasize that,
in particular for the lower tail,
this strategy does not require any control
of the variance of $W_n$.
\end{itremark}

\begin{pf*}{Proof of Lemmas~\ref{lem-AG-lower} and~\ref{lem-AG-upper}}
As in~\cite{AG} this is an application of Lemma~2.3 of~\cite{AG}.
For the lower tail, using the exponential Chebyshev's inequality,
the independence between $W_n$ and $L_n$, formula~(\ref{tool-lower})
and centering the random variables, we get
%
\begin{equation}
P(W_n < \xi_n) \leq \frac{E [e^{-\lambda(M_n - E[M_n])} ]} {
E [e^{-\lambda(L_n - E[L_n])} ]}
e^{-\lambda(E[M_n] - E[L_n] - \xi_n - c_n)}.
\end{equation}
With, for all $t\in{\mathbb R}$, $f(t) = e^t - (1+t)$
and $g(t) = (e^t - 1)^2$, by Lemma~2.3 of~\cite{AG},
%
\begin{eqnarray}
&&\frac{E [e^{-\lambda(M_n - E[M_n])} ]} {
E [e^{-\lambda(L_n - E[L_n])} ]}
\nonumber
\\[-8pt]
\\[-8pt]
\nonumber
&&\qquad \leq \exp \Biggl\{ f(-\lambda) \bigl(E[M_n] -
E[L_n]\bigr) + \frac{\kappa}{2} g(-\lambda) \sum
_{i=1}^{N_n}E^2\bigl[Y_i^{(n)}
\bigr] \Biggr\}.
\end{eqnarray}
We conclude by observing that for all $t\in{\mathbb R}$,
%
\begin{equation}
f(t) \leq\frac{t^2}{2}e^{[t]_+}\quad \mbox{and}\quad g(t)\leq
t^2e^{2[t]_+},
\end{equation}
where $[\cdot]_+$ stands for the positive part. The proof
for the upper tail is similar.
\end{pf*}

\subsection{On a discrete
mean value property of Green's function}

\begin{itproposition}\label{gauss}
Consider $d\geq2$.
There is a constant $K_d$
such that, for any~$n$ and $R$
with $n-n^{{1}/{3}}\leq R \leq n$
and $z$ in ${\mathbb B}(0,R)$
with $n - \|z\| \leq1$,
%
\begin{equation}
\biggl\llvert \bigl|{\mathbb B}(0,R)\bigr| G_n(0,z) - \sum
_{y\in{\mathbb B}(0,R)} G_n(y,z) \biggr\rrvert \leq
K_d.
\end{equation}
\end{itproposition}

\begin{pf}
For $n-R$ large enough
(larger than some constant
that depends only on $d$)
this is Theorem 5.2
of~\cite{AG}.
For $n= R$ this is a direct consequence
of Lemmas~2 and~3 of~\cite{lawler95}.
For the remaining cases, one can use
the same Lemmas in conjunction
with Lemma~5 of~\cite{lawler95}.
\end{pf}

\begin{itremark}\label{rem-2}
For the inner bound we will use Proposition~\ref{gauss}
with $R = n$.
For the outer bound we will use Proposition~\ref{gauss}
with $n-R$ of order $\log n$ in dimension 2
and $\sqrt{\log n}$ in dimension $d\geq3$.
\end{itremark}
%
\section{Inner error}\label{sec-inner}
\subsection{Exploration by waves}
We choose the following height sequence. For any positive integer $n$,
$h(n)=\sqrt{\log(n)}$ in $d\ge3$, and $h(n)=\log(n)$ in \mbox{$d=2$}.
We partition $\Z^d$ into concentric shells of heights $h(n)$.
We define $\S_0=\B(0,h(n))$, and for $k\ge1$,
\begin{equation}
\label{def-S} \S_k=\B\bigl(0,(k+1)h(n)\bigr)\bs\B\bigl(0,kh(n)\bigr)
\quad\mbox{and}\quad \Sigma_k=\p\B\bigl(0,kh(n)\bigr).
\end{equation}

We realize the internal DLA
with $N =|\B(0,n)|$ explorers
as an exploration wave process, where
concentric shells are covered in turn; see Section 3 of~\cite{lawler95}.

We fix an integer $k$. For a site $z\in\Sigma_{k}$,
we call \textit{cell} centered on $z$, $\C(z):=\B(z,h(n))\cap\S_{k}$,
and we call \textit{tile} centered on $z$, $\T(z):=\B(z,h(n)/2)\cap
\Sigma_k$.
A~generic cell is denoted $\C$, and a generic tile is denoted $\T$.
Note the obvious facts
\begin{equation}
\label{obvious-1} \bigcup_{z\in\Sigma_k} \B\bigl(z,h(n)\bigr)
\supset\S_k.
\end{equation}
Before covering shell $\S_k$, one
stops the unsettled explorers on $\Sigma_k$.
Following~\cite{AG},
for $z\in\Sigma_k$,
we prove that the $W_{kh(n)}(N\ind_0,\T)$
explorers stopped on $\T= \T(z)$
are likely to cover $\C(z)$,
if
$kh(n)\leq n-Ah(n)$ for a large enough constant $A$.
More precisely, we show that the probability
of the event $\{\S_k\not\subset A(N)\}$
is smaller, for $A$ large enough,
than any given power of $1/n$.
As first observed in \cite{lawler92},
\begin{equation}
\label{lawler-1} W_{kh(n)}(N\ind_0,\T)+M_{kh(n)}\bigl(
\B\bigl(0,kh(n)\bigr),\T\bigr) \ge\tilde M_{kh(n)}(N\ind_0,
\T).
\end{equation}
Since (\ref{lawler-1}) corresponds
to an inequality of type (\ref{tool-lower}),
we wish to use Lem\-ma~\ref{lem-AG-lower},
but we need to ensure (H1) and (H2).

First, if $\tilde\B(r)$ denotes the sites of $\B(0,kh(n))$ at a
distance less than $r$
from $\T$, there is $L$ and $\rho_d>1$ (which depend only on the
dimension),
such that
\begin{equation}
\label{checkH1} \sup_{y\in\B(0,kh(n))\bs\tilde\B(Lh(n))} P_y\bigl(S\bigl(H(
\Sigma_k)\bigr)\in\T\bigr) < \frac{\rho_d-1}{\rho_d};
\end{equation}
(see Lemma 5.1 of \cite{AG}).
Set $c_n=|\tilde\B(Lh(n))|$, and note that $c_n\le c(L h(n))^d$
for some constant $c$.
From (\ref{lawler-1}) we have
\begin{eqnarray}
\label{lawler-3} &&W_{kh(n)}(N\ind_0,\T)+ M_{kh(n)}
\bigl(\B\bigl(0,kh(n)\bigr)\bs\tilde\B\bigl(Lh(n)\bigr),\T\bigr)
\nonumber
\\[-8pt]
\\[-8pt]
\nonumber
&&\qquad\ge \tilde
M_{kh(n)}(N\ind_0,\T)- c_n.
\end{eqnarray}
We will use Lemma~\ref{lem-AG-lower}
with $L_n = M_{kh(n)}(\B(0,kh(n))\bs\tilde\B(Lh(n)),\T)$
and we note that (H1) is ensured by~(\ref{checkH1}).
Let us define
\begin{equation}
\label{def-mu}\qquad \mu(\T)=E \bigl[ {M_{kh(n)}(N\ind_0,\T)}
\bigr]-E \bigl[ {M_{kh(n)}\bigl(\B\bigl(0,kh(n)\bigr)\bs\tilde
\B_\T\bigl(Lh(n)\bigr),\T\bigr)} \bigr].
\end{equation}
We consider the event that $\S_k$ is not covered, and use the
bound
%
\begin{eqnarray}
\label{goal-inner}
&&P ( {\S_k\mbox{ not covered}} )\nonumber\\
&&\qquad\le P \bigl( {
\exists\T\subset\Sigma_k \dvtx W_{kh(n)}(N\ind_0,\T)<
\tfrac{1}{3} \mu(\T)} \bigr)
\\
&&\qquad\quad{}+ P \bigl( {\S_k\mbox{ not covered}, \forall\T\subset
\Sigma_k\dvtx W_{kh(n)}(N\ind_0,\T)\ge
\tfrac{1}{3} \mu(\T)} \bigr).\nonumber
\end{eqnarray}
In the next sections, we compute $\mu(\T)$, and estimate the probabilities
of the two events on the right-hand side of (\ref{goal-inner}).

\subsubsection{Mean number of explorers crossing a tile}
If $\mathcal{T}$ is a tile of a cell $\mathcal{C}$ which belongs to
shell $\S_k\subset\B(0,n)$, at a distance $A h(n)$ from $\B
(0,n)$, then
we show that for some positive constants $\{c_d, d\ge2\}$,
\begin{equation}
\label{main-1} \mu(\T)\ge c_d A h(n)^d.
\end{equation}
The inequality in (\ref{main-1})
follows as in \cite{AG}, Section 4.2,
and relies on Proposition~\ref{gauss}.
Note that~(\ref{main-1}) ensures (H2).

\subsubsection{$W_{kh(n)}(N\ind_0,\T)$ is unlikely to be small}
Like in (4.17) and (4.18) of Section 4.2 of \cite{AG},
there are constants $C_d$ such that
\begin{equation}
\label{cond-H2}
\sum_{y\in\B(0,kh(n))} P_y^2
\bigl( {S\bigl(H(\Sigma_k)\bigr)\in\T} \bigr)\le \cases{
C_2 h^2(n)\log(n), & \quad$\mbox{for } d=2,$
\vspace*{2pt}\cr
C_d h^d(n), & \quad$\mbox{for } d\ge3.$ }
\end{equation}
By Lemma~\ref{lem-AG-lower},
since for $A$ large enough
we have $\mu(\T)\geq3cL^dh(n)^d\geq3c_n$,
there are
positive constants $\{c'_d, d\ge2\}$ such that
\begin{eqnarray}
\label{outer-1}
&& P \bigl( {W_{kh(n)}(N\ind_0,\T)<
\tfrac{1}{3} \mu(\T)} \bigr)
\nonumber
\\[-8pt]
\\[-8pt]
\nonumber
&&\qquad\le \cases{
\exp \bigl( {-\lambda\kappa_2 A h^2(n)+
\lambda^2 c'_2 h^2(n)\log (n)}
\bigr),& \quad$\mbox{for } d=2,$
\vspace*{2pt}\cr
\exp \bigl( {-\lambda\kappa_d A h^d(n)+
\lambda^2 c'_d h^d(n)} \bigr),&\quad
$\mbox{for } d\ge3.$ }
\end{eqnarray}
Thus, after optimizing over $\lambda$, we get
\begin{eqnarray}
\label{outer-11}
&& P \biggl( {\exists z\in\Sigma_k\dvtx  W_{kh(n)}
\bigl(N\ind_0,\T(z)\bigr)<\frac
{1}{3} \mu(\T)} \biggr)
\nonumber
\\[-8pt]
\\[-8pt]
\nonumber
&&\qquad\le
\cases{ %
n^2\exp \biggl( {-
\displaystyle\frac{\kappa_2^2 A^2 h^2(n)}{4c'_2 \log(n)}} \biggr),&\quad $\mbox{for } d=2,$
\vspace*{2pt}\cr
n^d\exp \biggl( {-\displaystyle\frac{\kappa_d^2 A^2 h^d(n)}{4c'_d}} \biggr),& \quad$\mbox{for } d\ge3,$}
\end{eqnarray}
and the event
$\{W_{kh(n)}(N\ind_0,\T)\le\frac{1}{3} \mu(\T)\}$ has a probability
that decreases, for $A$ large enough, faster than any given power of $1/n$.

\subsubsection{$\C$ is likely to be covered when $W_{kh(n)}(N\ind_0,\T)$ is large}
We consider here the event $\{\forall\T\subset\S, W_{kh(n)}(N\ind_0,\T)\ge
\frac{\kappa}{3} Ah^d(n)\}$.
Consider shell $\S_k$ at a distance $Ah(n)$ from $\p\B(0,n)$.
Since $\S_k$ is the union of $\B(z,h(n))$ when $z\in\Sigma_k$,
Lemma~\ref{lem-ld} implies, when $d=2$, that
\begin{eqnarray}
\label{main-d2}
&& P \biggl(\S_k\notin A(N) \mbox{ and }
W_{kh(n)}(N\ind_0,\T) >\frac{1}{3}\mu(\T) \mbox{ for
all }\T \biggr)
\nonumber
\\[-8pt]
\\[-8pt]
\nonumber
&&\qquad\le|\S_k| \exp \biggl( {-\kappa_2 \kappa A
\frac{h^2(n)}{\log(n)}} \biggr).
\end{eqnarray}
%
We obtain a bound smaller than any power of $1/n$ when
$h(n)=\log(n)$
and $A$ is large enough. When $d\ge3$, then we have
\begin{eqnarray}
\label{main-d3} &&P \bigl(\S_k\notin A(N)\mbox{ and
}W_{kh(n)}(N\ind_0,\T)>\tfrac{1}{3}\mu(\T) \mbox{ for
all }\T \bigr)
\nonumber
\\[-8pt]
\\[-8pt]
\nonumber
&&\qquad\le|\S_k| \exp \bigl( {-\kappa_d \kappa
Ah^2(n)} \bigr).
\end{eqnarray}
For any given power of $n$, we obtain a negligible bound
when $h^2(n)=\log(n)$
and~$A$ is large enough.

\section{Outer error}\label{sec-outer}
In this section,
we prove the outer error estimate~(\ref{theo-res2}).
This is a consequence of our inner error estimates,
of Lemma~\ref{lem-levine}, combined with
coupling with a flashing process of \cite{AG}. When dimension $d=2$,
and for $A$ large to be chosen later,
we decompose the event $\{\delta_O(n)\ge A\log(n)\}$, as
\begin{equation}
\label{decomp-outer}\quad \bigl\{ {\delta_O(n)\ge A\log(n)} \bigr\}=
\bigcup_{i\ge1} \bigl\{ {\delta_O(n)\in\bigl[A
\log(n)+i-1,A\log(n)+i\bigr[} \bigr\}.
\end{equation}
In dimension $d\ge3$, $\sqrt{\log(n)}$ replaces $\log(n)$ in (\ref
{decomp-outer}).
Note that the index $i$ is at most of order $n^d$.
Now, we fix $i\ge1$, and we set
$3h(n)=A\sqrt{\log(n)}+i$ in $d\ge3$, and $3h(n)=A\log(n)+i$ in $d=2$.
We now consider the event $\{\delta_O(n)\in[3h(n)-1,3h(n)[\}$. We also define
\[
\Sigma=\B\bigl(0,n+3h(n)\bigr)\bs\B\bigl(0,n+3h(n)-1\bigr).
\]
Note now that
\begin{eqnarray}
\label{start-1}
&& P \bigl( {\delta_O(n)\in\bigl[3h(n)-1,3h(n)\bigr[} \bigr)
\nonumber
\\[-8pt]
\\[-8pt]
\nonumber
&&\qquad\le P
\biggl( {\bigcup_{ z\in\Sigma} \bigl\{ {z\in A(N), \delta_O(n)=\| z\|-n} \bigr\}} \biggr).
\end{eqnarray}
For $z\in\Sigma$, and in view of Lemma~\ref{lem-levine}, we define
\begin{equation}
\label{def-G}\qquad G(z)= \bigl\{ {z\in A(N), \delta_O(n)=\|z\|-n, \bigl|A(N)
\cap\B \bigl(z,h(n)\bigr)\bigr|>\beta h^d(n)} \bigr\}.
\end{equation}
To prove that $P(z\in A(N),\delta_O(n) =\|z\| - n)$ is smaller
than any given power of~$1/n$, we
further split the event
into two pieces:
\begin{eqnarray}
\label{start-2}
&&P \bigl( {z\in A(N), \delta_O(n)=\|z\|-n} \bigr)
\nonumber
\\[-8pt]
\\[-8pt]
\nonumber
&&\qquad\le
P\bigl(G(z)\bigr)+ P \bigl( {z\in A(N), \bigl|A(N)\cap\B\bigl(z,h(n)\bigr)\bigr|\le\beta
h^d(n)} \bigr).
\end{eqnarray}
The second term on the right-hand side of (\ref{start-2}) is dealt
with using
Lemma~\ref{lem-levine}.
We deal now with $G(z)$. Note that under $\{ \delta_O(n)\in
[3h(n)-1,3h(n)[\}$, no
explorer escapes $\B(0,n+3h(n))$. Thus, on $G(z)$, there are at least
$\beta h^d(n)$ explorers which settle on $\B(z,h(n))$ before exiting
$\B(0,n+3h(n))$.
We now express the event $G(z)$ in term of \textit{flashing explorers}, as
introduced in \cite{AG}.

\subsection{On a flashing process}
We refer the reader to Section 3.1 of \cite{AG} for a definition
of flashing processes.
Here, we partition $\Z^d$ into shells encaging $\B(0,n)$, with for
$k\ge0$,
\[
\S_{k}=\B\bigl(0,n+2(k+1)h(n)\bigr)\bs\B\bigl(0,n+2kh(n)\bigr).
\]
Also, for $k\ge0$, let $\Sigma_k=\p\B(0,n+(2k+1)h(n))$.
We now consider the flashing process.
Explorers behave like internal DLA explorers, as long as they
stay in $\B(0,n)$. After exiting $\B(0,n)$ they do not
flash until their hitting of $\Sigma_0$, and behave like
\textit{flashing explorers} as defined in Section 3.1 of \cite{AG}.
In shells $\{\S_k, k\ge0\}$, cells and
tiles have the meaning given in Section 4 of \cite{AG}.
The key features the reader has to keep in mind are as follows:
\begin{itemize}
\item If a flashing explorer is unsettled up to time $H( \Sigma_k)$,
then after time $H(\Sigma_k)$, it probes one site distributed
almost uniformly
over the \textit{cell} centered at $S(H(\Sigma_k))$, and
settles if the site is unoccupied.
\item When an explorer leaves the cell centered on $S(H(\Sigma_k))$,
it cannot afterward settle in $\S_k$, but perform a simple random walk,
independent of other explorers, until it hits $\Sigma_{k+1}$. Thus, if
we know that an explorer has reached at time $t$ a site of
$\B(0,n+(2k+1)h)\bs\B(0,n+2kh)$, then it performs after time $t$ a simple
random walk, independent of its surroundings, until it reaches $\Sigma_k$.
\item We can build the internal DLA cluster, $A(N)$,
and the flashing cluster $A^*(N)$ using the same trajectories
$S_1,\ldots,S_N$ such that
\begin{equation}
\label{feature-coupling} A(N)=\bigcup_{i=1}^N
\bigl\{S_i\bigl(T(i)\bigr)\bigr\}\quad \mbox{and}\quad A^*(N)=\bigcup
_{i=1}^N \bigl\{S_i
\bigl(T^*(i)\bigr)\bigr\},
\end{equation}
and for all $i=1,\ldots,N,  T^*(i)\ge T(i)$. This last property,
at the heart of our coupling argument between the flashing process
and the original DLA, is
fundamental. It implies that if a DLA explorer has crossed a site before
settling, then the corresponding flashing explorer has also crossed the
site before settling.
\end{itemize}

Before introducing more notation, let us explain the simple idea behind
our estimate.

\textit{Heuristics}.
Using representation (\ref{feature-coupling}), event $G(z)$ for $z\in
\Sigma$
implies that at least $\beta h^d(n)$
\textit{flashing explorers} hit $\B(z,h(n))$ before exiting $\B(0,n+3h(n))$.
Consider these explorers after the moment they enter $\B
(z,h(n))\subset\S_1$
for the first time. They
are behaving as independent random walks until they hit $\Sigma_1$.
Now, a fraction
must hit $\Sigma_1$ on $\B(z,2h(n))\cap\Sigma_1$. We show that this
latter event
has a probability we can estimate through the approach of \cite{AG}.

For simplicity, let us call $R_1$ the radius of $\Sigma_1$, that
is, $R_1=n+3h(n)$. Recall that for $\Lambda\subset
\B(0,R_1)\cup\p\B(0,R_1)$, we call
$W_{R_1}(N\ind_0,\Lambda)$ the number of flashing explorers which
hit $\Lambda$ before (or as)
they hit $\Sigma_1$. In this section, the initial configuration
is always $N\ind_0$, and we omit this coordinate in $W_{R_1}$
to simplify notation. Under our coupling
(\ref{feature-coupling}), we have
\begin{equation}
\label{outer-3} G(z)\subset \bigl\{ {W_{R_1} \bigl(\B\bigl(z,h(n)\bigr)
\bigr)\ge\beta h^d(n)} \bigr\}.
\end{equation}

Let $z'$ be the closest site of $\Sigma_1$
to the line $(0,z)$, and note that $\|z-z'\|\le1$.
Note that a fraction of the
$W_{R_1}(\B(z,h(n)))$ independent random walks in $\B(z,h(n))\cap\B
(0,n+3h(n))$,
must hit $\Sigma_1$ in a neighborhood of $z'$.
Indeed, first note that since $z'\in\Sigma_1$, we have
\begin{equation}
\label{def-L}
\bigl|\p\B\bigl(z',2h(n)\bigr)\cap\B\bigl(0,n+3h(n)\bigr)\bigr|
\ge\tfrac{1}{4} \bigl|\p\B\bigl(z',2h(n)\bigr)\bigr|.
\end{equation}
Now, for any $y\in\p\B(z,h(n))$, a random walk starting on $y$,
exits $\B(z',2h(n))$
on any site of $\p\B(z',2h(n))$ with a probability proportional to
$(2h(n))^{1-d}$.
Thus, there is a positive constant $\rho$ such that
\begin{equation}
\label{exit-uniform} \inf_{y\in\p B(z,h(n))} P_y \bigl( {S \bigl( {H
\bigl(\p\B\bigl(z',2h(n)\bigr)\bigr)} \bigr)\in\S_2}
\bigr)\ge\rho.
\end{equation}
In other words, each flashing explorer stopped on $\p\B(z,h(n))$
before hitting $\Sigma_1$
has a probability at least $\rho$ to exit $\Sigma_1$
from $\Sigma_1\cap\B(z,2h(n))$. Thus, there is a positive constant
$I$, such that for any large enough integer $k$,
\begin{equation}
\label{outer-4} \qquad P \biggl( {W_{R_1} \bigl( {\B\bigl(z',2h(n)
\bigr)\cap\Sigma_1} \bigr)<\frac
{\rho}{2} k \Big| W_{R_1}
\bigl(\B\bigl(z,h(n)\bigr)\bigr)>k} \biggr)\le\exp ( {- I k} ).
\end{equation}
From (\ref{outer-3}), we have
%
\begin{eqnarray}
\label{outer-5}
\bigcup_{z\in\Sigma} G(z)&\subset& \bigcup
_{z'\in\Sigma_1} \biggl\{ {W_{R_1} \bigl( {\B
\bigl(z',2h(n)\bigr)\cap\Sigma_1} \bigr)\ge
\frac{\rho}{2} \beta h^d(n)} \biggr\}
\nonumber
\\
&&{}\cup\bigcup_{z'\in\Sigma_1} \biggl\{ W_{R_1} \bigl(
{\B\bigl(z',2h(n)\bigr)\cap\Sigma_1} \bigr)<
\frac
{\rho}{2} \beta h^d(n) \mbox{ and}\\
&&\hspace*{79pt}\qquad W_{R_1} \bigl( {
\B\bigl(z,h(n)\bigr)} \bigr)\ge\beta h^d(n) \biggr\}.\nonumber
\end{eqnarray}
Let us now define, for any $a>0$,
\begin{equation}
\label{def-F} F(a)= \bigcup_{z\in\Sigma_1} \bigl\{
{W_{R_1} \bigl( {\B\bigl(z,2h(n)\bigr)\cap \Sigma_1} \bigr)
\ge ah^d(n)} \bigr\}.
\end{equation}
Thus, from (\ref{outer-5}) and (\ref{outer-4}), and for some constant $C>0$
%
\begin{eqnarray}
\label{outer-6} P \biggl( {\bigcup_{z\in\Sigma} G(z) }
\biggr)&\le &P \biggl( {\bigcup_{z'\in\Sigma_1} \biggl\{
{W_{R_1} \bigl( {\B \bigl(z',2h(n)\bigr)\cap
\Sigma_1 } \bigr)\ge\frac{\rho}{2} \beta h^d(n) }
\biggr\} } \biggr)
\nonumber\\
&&{}+ |\Sigma_1|\sup_{z'\in\Sigma_1} P \biggl( W_{R_1} \bigl(
{\B\bigl(z',2h(n)\bigr)\cap\Sigma_1 } \bigr)<
\frac
{\rho}{2} \beta h^d(n)\Big|
\nonumber
\\[-8pt]
\\[-8pt]
\nonumber
&&\hspace*{95pt}\qquad  W_{R_1} \bigl( {\B
\bigl(z,h(n)\bigr)} \bigr)\ge \beta h^d(n)  \biggr)
\\
&\le &P \biggl( {F \biggl( {\frac{\rho}{2} \beta } \biggr) } \biggr)+C
n^{d-1} \exp \biggl( {- I \frac{\rho}{2} \beta h^d(n)}
\biggr).\nonumber
\end{eqnarray}
It remains to show that for any fixed $a$, we can find $A$ [defining
$h(n)$] such that
$P(F(a))$ is smaller than any given power of $1/n$.

\subsection{Estimating $P(F(a))$}
Note that by definition of $\delta_I(n)$,
for $z\in\Sigma_1
=\p\B(0,n+3h(n))$ and
$\T_z = \B(z, 2h) \cap\Sigma_1$,
$W_{R_1}(\T_z)$ satisfies the inequality
\begin{equation}
\label{ineq-w1} W_{R_1}(\T_z)+M_{R_1}\bigl(\B
\bigl(0,n-\delta_I(n)\bigr),\T_z\bigr)\le \tilde
M_{R_1}(N\ind_0,\T_z).
\end{equation}
Thus, for some large constant $\alpha_d$ to be chosen later, we have
\begin{eqnarray}
\label{ineq-alex} &&\ind_{\delta_I(n)\le\alpha_d {h(n)}/{(2A)}} \biggl( {W_{R_1}(
\T_z)+M_{R_1}\biggl(\B\biggl(0,n-\alpha_d
\frac{h(n)}{2A}\biggr)\biggr)} \biggr)
\nonumber
\\[-8pt]
\\[-8pt]
\nonumber
&&\qquad\le\tilde M_{R_1}(N
\ind_0,\T_z).
\end{eqnarray}

Inequality (\ref{ineq-w1}) puts us
in the setting of Lemma~\ref{lem-AG-upper}. Thus, we first need
to compute
\begin{equation}
\label{tildemu-0}\qquad \tilde\mu(\T_z)=E \bigl[ {M_{R_1}(N
\ind_0,\T_z)} \bigr] -E \biggl[ {M_{R_1}\biggl(
\B\biggl(0,n-\alpha_d \frac{h(n)}{2A}\biggr),\T_z
\biggr)} \biggr].
\end{equation}
Following the same computations
as in Section 4.3 of \cite{AG}, we have for some constant $K$
\begin{equation}
\label{tildemu-1} \tilde\mu(\T_z)\le K \biggl( {\alpha_d
\frac{h(n)}{A}n^{d-1}} \biggr)\times\frac{h^{d-1}(n)} {
n^{d-1}}\le
\frac{K\alpha_d}{A}h^d(n).
\end{equation}
%
Second, note that as in Section 4.3
of \cite{AG}, we have that for constants \mbox{$\{c_d, d\ge2\}$},
\begin{equation}
\label{tildemu-3} \sum_{z\in\B(0,n)} P_z^2
\bigl( {S\bigl(H(\Sigma_1)\bigr)\in\T_z} \bigr)\le
\cases{
 c_2 h^2(n)\log(n),
&\quad$\mbox{if } d=2,$
\vspace*{2pt}\cr
c_d h^d(n), &\quad$\mbox{if } d\geq3.$ }
\end{equation}
In optimizing over $\lambda$ in (\ref{lem-AG-upper}), we find
for (other) constants $\{c_d, d\ge2\}$, if $A$ is choosen large enough
\begin{eqnarray}
\label{outer-upper}\qquad && P \bigl( {\exists z\in\Sigma_1 \dvtx
W_{R_1}(\T_z)\ge a h^d(n)} \bigr)
\nonumber
\\[-8pt]
\\[-8pt]
\nonumber
&&\qquad\le P \biggl( {
\delta_I(n)> \alpha_d \frac{h(n)}{2A}}
\biggr)+n^{d} \cases{ %
\exp \biggl(
{-c_2 \displaystyle\frac{h^2(n)}{\log(n)}} \biggr), &$\mbox{if } d=2,$
\vspace*{2pt}\cr
\exp \bigl( {-c_d h^d(n)} \bigr), &$\mbox{if } d\geq3.$}
\end{eqnarray}
We conclude using the fact,
for $\alpha_d$ large enough,
the first term of the sum
in the right-hand side of~(\ref{outer-upper})
is smaller
than any given power of $1/n$.
This was proved in Section~\ref{sec-inner}
for the original internal DLA and the same proof
can be adapted for the flashing process we consider
here.
The only difference is that we need
a stronger version of Lemma~\ref{lem-ld}
where $P(\B(0,R) \not\subset A(\eta))$
is replaced by $P(\B(0,R) \not\subset A_{\alpha R}(\eta))$
for some large $\alpha$
(this stronger version of the lemma
is actually what we prove in the \hyperref[app]{Appendix}).
Indeed, we can use Lemma~\ref{lem-ld}
in the context of our flashing process
by considering only explorers
that do not exit $\B(0,n)$.
Once $\alpha_d$ is fixed, we choose $A$ large enough so that (\ref{outer-upper}) holds.

\begin{appendix}\label{app}

\section{\texorpdfstring{Proof of Lemma~\lowercase{\protect\ref{lem-ld}}}{Proof of Lemma 1.3}}

We fix $\eta$, a configuration of $AR^d$ explorers
in $\B(0,R/2)$, and we choose $z\in\B(0,R)$. Then
\begin{eqnarray}
\label{idea-lawler} P \bigl( {\B(0,R)\not\subset A(\eta)} \bigr)&\le& P \bigl( {
\B(0,R)\not\subset A_{\alpha R}(\eta)} \bigr)
\nonumber
\\[-8pt]
\\[-8pt]
\nonumber
&\le& \sum
_{z\in\B(0,R)} P \bigl( {W_{\alpha R}(\eta,z)=0} \bigr)
\end{eqnarray}
for any $\alpha> 1$ (in the sequel $\alpha$
will have to be taken large enough). Let $L$ be a large positive
real to be fixed later, and let
$\zeta$ be the configuration with one explorer on each
site of $\B(0,\alpha R)\bs\B(z,L)$. We have
\begin{equation}
\label{lawler-2} W_{\alpha R}(\eta,z)+M_{\alpha R}(\zeta,z)\ge\tilde
M_{\alpha
R}(\eta,z) -\bigl|\B(z,L)\bigr|.
\end{equation}
Note that $W_{\alpha R}(\eta,z)$ and
$M_{\alpha R}(\zeta,z)$ are independent: we are in
the setting of Lemma~\ref{lem-AG-lower}.
Assume for a moment that conditions (H1) and (H2) hold, and in addition,
\begin{eqnarray}
\label{hyp-AG}&& E \bigl[ {M_{\alpha R}(\eta,z)} \bigr]-E \bigl[
{M_{\alpha R}(\zeta,z)} \bigr]
\nonumber
\\[-8pt]
\\[-8pt]
\nonumber
&&\qquad\ge \max \biggl( {3\bigl|\B(z,L)\bigr|,\sum
_{y\in\B(0,\alpha R)} P_y ( {H_z< H_{\alpha R}}
)^2} \biggr).
\end{eqnarray}
Then, we have
\begin{equation}
\label{AG-main} P \bigl( {W_{\alpha R}(\eta,z)=0} \bigr) \le\exp \bigl( {-C
\bigl( {E \bigl[ {M_{\alpha R}(\eta,z)} \bigr]-E \bigl[ {M_{\alpha R}(
\zeta,z)} \bigr]} \bigr)} \bigr).
\end{equation}
We next consider separately the case $d\ge3$ and the case $d=2$,
estimate the expectation of
$\tilde M_{\alpha R}(\eta,z)-M_{\alpha R}(\zeta,z)$ and show (\ref{hyp-AG}).

\subsection{\texorpdfstring{The case $d\ge3$}{The case d >= 3}}
We show in this section that for some $\kappa_d>0$, and $A$ large enough,
\begin{equation}
\label{mean-5} E\bigl[\tilde M_{\alpha R}(\eta,z)-M_{\alpha R}(\zeta,z)
\bigr] \ge\frac{\kappa_d}{2} A R^2\gg3\bigl|\B(z,L)\bigr|.
\end{equation}
The proof is based on the following classical estimates.
There are $a_1,a_2$ positive constants such that
for any $y,z\in\Z^d$
\begin{equation}
\label{hitting-estimate} \frac{a_1}{1+\|y-z\|^{d-2}}\le P_y(H_z<
\infty)\le \frac{a_2}{1+\|y-z\|^{d-2}}.\vadjust{\goodbreak}
\end{equation}
Note first that when $L$ is large enough, (H1) holds. Indeed,
\begin{eqnarray}
\label{hypo-folk} \sup_{y: \|z-y\|>L} P_y(H_z<H_{\alpha R})
&\le& \sup_{y: \|z-y\|>L} P_y(H_z<\infty)
\nonumber
\\[-8pt]
\\[-8pt]
\nonumber
&\le&
\frac{a_2}{1+L^{d-2}}\le\frac{\kappa-1}{\kappa} \quad\mbox{with } \kappa>1.
\end{eqnarray}
We now estimate the mean number of explorers hitting $z$.
%
\begin{eqnarray}
\label{mean-1}
&& E\bigl[M_{\alpha R}(\eta,z)\bigr]-E\bigl[M_{\alpha R}(
\zeta,z)\bigr]\nonumber\\
&&\qquad = \sum_{y\in\B(0,R/2)}\eta(y) P_y (
{H_z<H_{\alpha R}} )-\sum_{y\in\B(0,\alpha R)\bs\B(z,L)}
P_y ( {H_z<H_{\alpha R}} )
\\
&&\qquad\ge \sum_{y\in\B(0,R/2)}\eta(y) P_y (
{H_z<H_{\alpha R}} )-\sum_{y\in\B(0,\alpha R)}P_y
( {H_z<\infty} ).\nonumber
\end{eqnarray}
Note that for $y\in\B(0,R/2)$, we have
%
\begin{eqnarray}
\label{mean-2} P_y ( {H_z<H_{\alpha R}} )&=&
P_y(H_z<\infty)-E_y \bigl[ {
\ind_{H_{\alpha R}<H_z} P_{S(H_{\alpha R})} ( {H_z<\infty} )} \bigr]
\nonumber\\
&\ge& \frac{a_1}{1+\|y-z\|^{d-2}}-E_y \biggl[ {\frac{a_2} {1+\|
S(H_{\alpha R})-z\|^{d-2}}} \biggr]
\\
&\ge& \inf_{y\in\B(0,R/2)} \frac{a_1}{1+\|y-z\|^{d-2}}- \sup_{y\in\p\B(0,\alpha R)}
\frac{a_2}{1+\|y-z\|^{d-2}}.\nonumber
\end{eqnarray}
Now, for a constant $\alpha$ which depends only on $a_1,a_2$,
there is $\kappa>0$ such that
\begin{equation}
\label{mean-3} \inf_{y\in\B(0,R/2)} P_y ( {H_z<H_{\alpha R}}
)\ge\frac
{\kappa}{R^{d-2}}.
\end{equation}
Now, using (\ref{mean-3}) in (\ref{mean-1}), we have
a constant $c$ such that
%
\begin{eqnarray}
\label{mean-4} E\bigl[\tilde M_{\alpha R}(\eta,z)-M_{\alpha R}(\zeta,z)
\bigr]&\ge& A R^d \frac{\kappa}{R^{d-2}}-\sum_{y:\|y-z\|<\alpha R}
\frac{a_1}{1+\|y-z\|^{d-2}}
\nonumber
\\[-8pt]
\\[-8pt]
\nonumber
&\ge& \kappa A R^2-c a_2 (\alpha R)^2.
\end{eqnarray}
When $A$ is chosen large enough, we obtain (\ref{mean-5}).

Finally, there are constants $\{C_d, d\ge3\}$ such
that for any $z\in\B(0,R)$
\begin{equation}
\label{FMW} \sum_{y\in\B(0,\alpha R)} P_y (
{H_z< H_{\alpha R}} )^2\le \cases{
C_3\alpha R, & \quad$\mbox{for } d=3,$
\vspace*{2pt}\cr
C_4 \log(\alpha R), &\quad $\mbox{for } d=4,$
\vspace*{2pt}\cr
C_d, &\quad $\mbox{for } d\ge5.$}
\end{equation}
Thus, hypothesis (\ref{hyp-AG}) holds.
\subsection{The case $d=2$}
We still have
\begin{eqnarray}
\label{hitting-d2} P_y(H_z<H_{\alpha R})&=&
\frac{G_{\alpha R}(y,z)}{G_{\alpha R}(z,z)}= \frac{G_{\alpha R}(z,y)}{G_{\alpha R}(z,z)}\quad \mbox{and}
\nonumber
\\[-8pt]
\\[-8pt]
\nonumber
 G_{\alpha R}(z,y)&=&E_z
\bigl[ {a\bigl(S(H_{\alpha R}),y\bigr)} \bigr]-a(z,y),
\end{eqnarray}
where the \textit{potential kernel} $a(.,.)$ replaces
Green's function.
Note that for $0\le\|z\|+ R<\alpha R$, we have two positive
constants $K_2$ and $K_2'$ such that
\begin{equation}
\label{loop}  K_2'\log(2\alpha R) \geq G_{B(z,2\alpha R)}
\geq G_{\alpha R}(z,z) \geq G_{B(z,R)}(z,z) \geq K_2\log(R),\hspace*{-35pt}
\end{equation}
by Proposition 1.6.6 of Lawler \cite{lawler}.
To estimate $G_{\alpha R}(z,y)$, we use
Theorem~4.4.4 of \cite{lawler-limic} which establishes that
for $z\not=0$ (with $\gamma$ the Euler constant),
\begin{equation}
\label{main-potential} \biggl|a(0,z) - \frac{2}{\pi}\log \bigl( {\|z\|} \bigr)-
\frac{2\gamma+\log(8)}{\pi} \biggr|\le\frac{K_g}{\|z\|^2}.
\end{equation}
Thus, for $y\in\B(0,\alpha R)$, $0\le\|z\|\le R$, and $y\not= z$
%
\begin{equation}
\label{star} \biggl\llvert G_{\alpha R}(z,y) - \frac{2}{\pi} E \biggl[
\log \biggl( \frac{\|S(H_{\alpha R})-z\|}{\|y-z\|} \biggr) \biggr] \biggr\rrvert
\leq2K_g.
\end{equation}
When $y\in B(0,R/2)$, we get
%
\begin{equation}
G_{\alpha R}(z,y) \geq\frac{2}{\pi} \log \bigl( {2(\alpha-1)/3}
\bigr)-2K_g.
\end{equation}
We choose $\alpha$ large enough so that for some constant $C_1$,
we have, for all $y$ in $B(0,R/2)$,
\begin{equation}
\label{green-1} G_{\alpha R}(z,y)\ge C_1.
\end{equation}
Formulas~(\ref{hitting-d2}), (\ref{loop}) and~(\ref{green-1}) together
imply that
\begin{eqnarray}
\label{expectation-1} E\bigl[M_{\alpha R}(\eta,z)\bigr]&=&\sum
_{y\in\B(0,R/2)} \eta (y)P_y(H_z<H_{\alpha R})\nonumber\\
&\ge& \frac{C_1}{K_2'\log(2\alpha R)}\sum_{y\in\B(0,R/2)} \eta(y)
\\
&=&
\frac
{C_1A R^2}{K_2'\log(2\alpha R)}.\nonumber
\end{eqnarray}
Using Lemma 3 of \cite{lawler95}, we have, for some positive constant $C_2$,
\begin{equation}
\label{expectation-2} E\bigl[M_{\alpha R}(\zeta,z)\bigr] \leq E
\bigl[M_{\alpha R}\bigl({\mathbb B}(0,\alpha R),z\bigr)\bigr] \leq
\frac{C_2 (\alpha R)^2}{\log(R)}.
\end{equation}

We need now to choose $L$ to have (H1) satisfied.
Note that for $y\not= z$, (\ref{star}) and (\ref{loop}) yields
\begin{equation}
\label{green-2}\qquad P_y(H_z<H_{\alpha R})\le
\frac{1}{K_2 \log(R)} E \biggl[ {\frac{2}{\pi}\log \biggl( {\frac{\|S(H_{\alpha R})-z\|}{\|
y-z\|}}
\biggr)+2K_g} \biggr].
\end{equation}
If $\|z-y\|>R/\log(R)$, we obtain, for some constant $C_3$,
\begin{equation}
\label{green-3} P_y(H_z<H_{\alpha R})\le
\frac{C_3 \log ( {(\alpha+1) \log
(R)}  )}{\log(R)}.
\end{equation}
When $R$ is large enough, we have that (H1) holds
for $L=R/\log(R)$. Note that $|\B(0,L)|$ is of order $R^2/\log(R)^2$
and is much smaller than $R^2/\log(R)$.

Finally we need to control the sum of second moments.
Simply note that, from~(\ref{expectation-2}),
\begin{equation}
\label{check-d2} \sum_{y\in\B(0,\alpha R)\bs\B(0,L)} P_y^2(H_z<H_{\alpha R})
\le E\bigl[M_{\alpha R}(\zeta,z)\bigr]\le\frac{C_2\alpha^2 R^2}{\log(R)}.
\end{equation}

\section{\texorpdfstring{Proof of Lemma~\lowercase{\protect\ref{lem-ag}}}{Proof of Lemma 1.6}}
We will choose an $h$ such that $R/2h$ is a positive integer.
We divide $\S=B(0,2R)\setminus B(0,R)$
into $R/2h$ concentric shells of height $2h$.
For $k=1,\ldots,R/2h$, define
%
\begin{eqnarray}
\S_k&\hspace*{1pt}=&\B\bigl(0,2R-2(k-1)h\bigr)\bs\B(0,2R-2kh)\quad \mbox{and}
\nonumber
\\[-8pt]
\\[-8pt]
\nonumber
\Sigma_k &:=& \partial\B\bigl(0,2R-(2k-1) h\bigr).
\end{eqnarray}
Also, we set $\S_0=\B(0,2R)^c$. Then, we start on $z\in\p\B(0,2R)$
a flashing explorer
associated with this partition with an explored region $V$. The flashing
setting is much simpler than the one introduced in Section 3.1 of \cite{AG}.
There is an underlying simple random walk, say $S^*$, and each
shell $\S_1,\S_2,\ldots$ is associated with a flashing site.
These flashing sites,
say $\{Z_k, 0\leq k\leq2R/h\}$ are obtained as follows.
We set $Z_0 = z$, and for $k\geq1$
we draw a continuous random variable $R_k$ on $[0,h]$
with density in $r\in[0,h]\mapsto{dr^{d-1}}/{h^d}$:
the flashing site $Z_k$ is the exit site from $\B(S^*(H(\Sigma_k)),R_k)$
after time $H(\Sigma_k)$.
Then, the explorer settles on the first flashing site in
${\mathcal S}\setminus V$.
The purpose of the flashing construction is that:
(i) the flashing site is distributed almost uniformly
inside the ball $\B(S^*(H(\Sigma_k)),h)$; and
(ii) $P_z(H(\B(0,R)) < H(V^c))$ is bounded above
by the probability that the explorer crosses $\S$.

For a small $\beta$ to be chosen later, we say that $y\in\Sigma_k$
has a \textit{dense neighborhood} if $|\B(y, h) \cap V| \geq\beta h^d$, and
we call $D_k$ their set. There is
$\kappa>0$ such that knowing that $S^*$ has crossed $D_1,\ldots,D_{k-1}$:
\begin{itemize}
\item if $S^*(H(\Sigma_k))\notin D_k$, then the probability that
$S^*$ does not settle in $\S_k$ is smaller than $\kappa\beta$;
\item the probability that $S^*(H(\Sigma_k))\in D_k$ is smaller than
$\kappa|D_k|/h^{d-1}$ (see Lemma~5 of \cite{lawler92})
uniformly over the position of the previous flashing site (in $\S_{k-1}$
or, exceptionally, on the border of $S_{k-1}$).
\end{itemize}

Now, the flashing explorer has crossed the annulus $\S$ if $Z_k\in V$
for all $k\geq1$.
In other words,
\begin{equation}
\label{alex-0} \bigl\{ {H\bigl(\B(0,R)\bigr) < H\bigl(V^c\bigr)}
\bigr\}\subset\bigcap_{k=1}^{R/2h} \{
{Z_k\in V} \}.
\end{equation}
By successive conditioning, we obtain
\begin{equation}
\label{alex-1} P_z \Biggl( {\bigcap_{k=1}^{R/2h}
\{ {Z_k\in V} \}} \Biggr)\le\prod_{k=1}^{R/2h}
\biggl( {\kappa\beta+\frac{\kappa|D_k|}{h_k^{d-1}}} \biggr).
\end{equation}
By the arithmetic--geometric inequality and (\ref{alex-0}), we obtain
%
\begin{equation}
P_z\bigl(H\bigl(\B(0,R)\bigr) < H\bigl(V^c\bigr)\bigr)
\leq \Biggl(\kappa\beta+ \frac{\kappa}{R/2h} \sum_{k=1}^{R/2h}
\frac{|D_k|}{h^{d-1}} \Biggr)^{R/2h}.
\end{equation}
Note that each $y\in D_k$ satisfies $|\B(y,h)\cap V|\ge\beta h^d$,
but each
site in $\B(y,h)\cap V$ is in the neighborhood of at most $h^{d-1}$
sites of $D_k$.
Thus for some $\kappa'$,
%
\begin{equation}
\sum_{k=1}^{R/2h}\frac{\beta|D_k| h^d}{h^{d-1}} \leq
\kappa'|V|\qquad \mbox{i.e., } \frac{1}{R/2h}\sum
_{k=1}^{R/2h}\frac{|D_k|}{h^{d-1}} \leq\frac
{2\kappa'|V|}{\beta R h^{d-1}}.
\end{equation}
We choose now $\beta$ such that $4\kappa\beta<1$, and we choose the smallest
$h$ such that $R/2h$ is a positive integer and
\begin{equation}
\label{alex-3} h\geq\max \biggl\{h_0, \biggl(\frac{2\kappa'|V|}{\beta^2 R}
\biggr)^{{1}/{(d-1)}} \biggr\}.
\end{equation}
This adds a constraint on $|V|$,
\begin{equation}
\label{alex-4} |V|\le\frac{\beta^2}{2^d\kappa'} R^d.
\end{equation}
Instead of including (\ref{alex-4}) as a condition of our lemma, we
find it
more convenient to note that the probability we estimate is always less
than 1,
so that we deal with the case where (\ref{alex-4}) is violated with the
constant $a_d$
of (\ref{avoiding-traps}).

\section{\texorpdfstring{Proof of Lemma~\lowercase{\protect\ref{lem-levine}}}{Proof of Lemma 1.5}}
Recall that $a_d$ and $\kappa_d$ are the constants appearing in
Lemma~\ref{lem-ag}.
We define a positive constant
\begin{equation}
\label{def-gamma} \gamma=\max \biggl( {1, \biggl( {\frac{2a_d}{\kappa_d}}
\biggr)^{d-1}} \biggr).
\end{equation}
Choose now $\beta>0$ such that $4^d\beta\gamma\le1$ and $h_0=R/4\ge1$.
Note that
%
\begin{equation}
\label{sstar} \gamma|\eta|\le\gamma\beta R^d\le h_0^d.\vadjust{\goodbreak}
\end{equation}
We build now, by induction, a random subdivision of $\B(z,R)$ into
shells of heights $h_0,h_1,\ldots,$ in which, respectively,
$N_0,N_1,\ldots$ explorers
of $A(\eta)$ have settled. We emphasize that the randomness comes from
$A(\eta)$,
and that the event $\{0\in A(\eta)\}$ imposes to have $N_i\ge\lfloor
h_i \rfloor$,
for $i\ge0$.
Assume that $h_1,\ldots,h_k$ have been defined such that
\begin{equation}
\label{cond-induction} h_k\geq1 \quad\mbox{and}\quad \sum
_{i=1}^k h_i< \frac{R}{2}.
\end{equation}
We define $h_{k+1}^d=\gamma N_k\leq\gamma|\eta|$, and, by~(\ref{sstar})
we have $h_{k+1}\le h_0$.
Note also that $h_{k+1}\ge1$. Indeed, necessarily $N_k\ge\lfloor
h_k\rfloor$, so that
$h_{k+1}^d\ge\gamma\lfloor h_k\rfloor\ge\lfloor h_k\rfloor$.
Since $\min(h_1,\ldots,h_{k+1})\ge1$, the number of steps
before we violate (\ref{cond-induction}), say $L$,
is finite.
Obviously $L\le R$. Note that since $h_{L}\le h_0$,
\begin{equation}
\label{choice-2} \frac{R}{2} \leq\sum_{i=1}^L
h_i \leq h_L + \sum_{i=1}^{L-1}
h_i \leq\frac{R}{4}+\frac{R}{2}.
\end{equation}
Thus, we define
\begin{equation}
\label{choice-3} h_{L+1}=R- \Biggl(\sum_{i=0}^L
h_i \Biggr)\ge0.
\end{equation}
For any choice of integers $l,n_0,\ldots,n_l$, the event
$\{L=l,N_0=n_0,\ldots,N_L=n_l\}$ implies that
$n_1+\cdots+n_{l}$ explorers have crossed a shell $B(z,R)\bs B(z,R-h_0)$
by stepping on at most $n_0$ explorers settled in it, that
$n_2+\cdots+n_{L}$ explorers have crossed shell $B(z,R-h_0)\bs B(z,R-h_0-h_1)$
with $n_1$ explorers settled in it, and so on and so forth. Using
Lemma~\ref{lem-ag},
the fact that $n_i\le\beta R^d$, $l\le R$
and the notation $\delta=\frac{1}{d-1}$, we reach the following estimate:
%
\begin{eqnarray}
\label{upper-ld-1}
&&P \bigl( {0 \in A(\eta)} \bigr)
\nonumber\\
&&\qquad\le\mathop{\sum_{l\le R, n_0,n_1,\ldots,n_l\leq|\eta|}}_
{\forall i, n_i \geq\lfloor h_i\rfloor}
P(L=l, N_0=n_0,
\ldots,N_L=n_l)
\nonumber
\\[-8pt]
\\[-8pt]
\nonumber
&&\qquad\le R \bigl(\beta R^d\bigr)^{R+1} \mathop{\sup_
{ l\le R, n_0,n_1,\ldots,n_l\leq|\eta|}}_
{\forall i, n_i \geq\lfloor h_i\rfloor}
e^{a_d \sum_{i=1}^{L} i n_{i}}\\
&&\qquad\quad{}\times \exp \Biggl( {-\kappa_d \sum
_{i=1}^{L} n_{i} \biggl( { \biggl( {
\frac
{h_0^d}{n_0}} \biggr)^\delta+\cdots+ \biggl( {\frac
{h^d_{i-1}}{n_{i-1}}}
\biggr)^\delta} \biggr)} \Biggr).\nonumber
\end{eqnarray}
Now, note that by the arithmetic--geometric inequality, for $1\le i \leq l$
(and using $h_i \leq h_0$)
%
\begin{eqnarray}
\label{choice-4} \qquad\frac{1}{i} \biggl( { \biggl( {\frac{h_0^d}{n_0}}
\biggr)^\delta +\cdots+ \biggl( {\frac {h_{i-1}^d}{n_{i-1}}} \biggr)^\delta}
\biggr)& \ge& \biggl( {\frac{h_0^d}{n_0}\times\cdots\times\frac
{h_{i-1}^d}{n_{i-1}}}
\biggr)^{\delta/i}
\nonumber
\\[-8pt]
\\[-8pt]
\nonumber
&=& \biggl( {\frac{h_0^d}{n_{i-1}}\gamma^{i-1}} \biggr)^{\delta/i} =
\biggl( {\frac{h_0^d}{h_i^d}\gamma^{i}} \biggr)^{\delta/i} \ge
\frac{2a_d}{\kappa_d}.
\end{eqnarray}
Thus, from (\ref{upper-ld-1}) and (\ref{choice-4}), we have
%
\begin{eqnarray}
\label{upper-ld-2} P \bigl( {0\in A(\eta)} \bigr) & \le &R \bigl(\beta
R^d\bigr)^{R+1} \mathop{\max_{ l\le R, n_0,n_1,\ldots,n_l\leq|\eta|}}_
{\forall i, n_i \geq\lfloor h_i\rfloor}
 \exp \Biggl(
{-a_d \sum_{i=1}^{L} i n_{i}} \Biggr)
\nonumber
\\[-8pt]
\\[-8pt]
\nonumber
&\le& R \bigl(\beta R^d\bigr)^{R+1} \mathop{\max_{l\le R,
n_0,n_1,\ldots,n_l\leq|\eta|}}_{
\forall i, n_i \geq\lfloor h_i\rfloor
} \exp
\Biggl( {-\frac{a_d}{\gamma} \sum_{i=1}^{L-1}
i h^d_{i+1}} \Biggr).
\end{eqnarray}
Since $h_1\le R/4$, note that we have $h_2+\cdots+h_L\ge R/4$
by (\ref{choice-2}).
By H\"older's inequality, note that for constants $\{c_d, d\ge2\}$,
\begin{eqnarray}
\label{holder} \qquad\sum_{i=1}^{L-1} i
h^d_{i+1}&\ge&\frac{
( {\sum_{i=1}^{L-1} h_{i+1}}  )^d}{
( {\sum_{i=1}^{L-1}{1}/{i^{1/(d-1)}}}  )^{d-1}}
\nonumber
\\[-8pt]
\\[-8pt]
\nonumber
&\ge &\cases{
c_2 \displaystyle\frac{R^2}{\log(L)}\ge c_2
\frac{R^2}{\log(R)},&\quad $\mbox{for }d=2,$
\vspace*{2pt}\cr
c_d \displaystyle\frac{R^d}{L^{d-2}}\ge c_d R^2, &\quad $\mbox{for } d\ge3.$ }
\end{eqnarray}
This completes the proof.
\end{appendix}

%


\printaddresses

\end{document}